\documentclass[11pt]{article}

\usepackage{graphicx}
\usepackage{amsthm}

\newcommand{\dfn}[1]{\textbf{#1}}
\newtheorem{statement}{Statement}
\newtheorem*{theorem}{Theorem}

\begin{document}

\title{Topology Explains Why Automobile Sunshades Fold Oddly}
\author{Curtis Feist and Ramin Naimi}
\date{}

\maketitle

\section*{Introduction}

Most readers are probably familiar with Magic Shade ``automatic folding'' sunshades,
even if they don't recognize the brand name.
These shades for automobile windshields are roughly the shape of a circular disk when ``open'',
and then fold up into a coil of several smaller circles when ``closed'',
thus allowing the shades to be stored in a relatively small space.

\begin{figure}[ht]

 \centering
 \includegraphics[width=60mm]{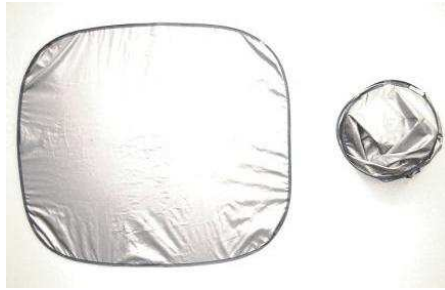}
 \caption{The shade in its open (left) and closed (right) positions.
 \label{photo-Figure}
 }
\end{figure}

Let's examine the structure of a Magic Shade a bit more closely.
The key part of the shade is a flat metal wire
that runs the length of the perimeter of the shade,
forming a continuous loop.
This wire is torsionally rigid.
That is, it doesn't mind being bent along its length,
but it does mind being twisted.
This observation will form the basis of our mathematical investigation,
as follows.
The wire frame of a Magic Shade,
when open, is in (roughly) the shape of a single circle,
with no twists in the wire.
When closed,
this wire frame must be coiled into several smaller circles/loops,
but \textit{without twisting the wire}.
Specifically, we ask the following questions:

\begin{enumerate}
\item
With the shade in its ``normal'' closed position,
how many loops does the wire make, and why?

\item
More generally,
what are \emph{all} possible closed positions for the shade,
in terms of how many loops the wire makes?

%Can the shade be put into other closed positions?
%If so, how many loops can the wire make, and why?
\end{enumerate}

Note that the ``how many'' part of the first question could be answered by observation,
but this would do nothing to address the ``why'' part of the question.

\begin{figure}[ht]

 \centering
 \includegraphics[width=20mm]{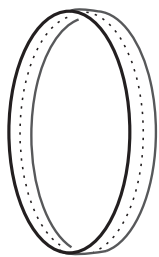}
 \caption{The band's core (dotted curve)
 and boundary circles (solid curves)
 with the shade in its open position.
 \label{open-band-Figure}
 }
\end{figure}

First, a bit of terminology: Let's call the flat wire frame a \dfn{band}.
This band has two edges,
which we will call \dfn{boundary circles}.
Let's call the circle midway between the boundary circles
(see Figure~\ref{open-band-Figure})
the \dfn{core} of the band.

Next, let's discuss what we mean by saying
that when the shade is in its closed position,
the band is ``coiled'' into several loops.
Informally, this would be similar to winding a length of ribbon
around a spool and then joining the ends together
(so that the piece of ribbon is itself a loop).
In fact, this is precisely what we have in mind,
except that for full generality
we also need to allow threading of the ribbon
beneath loops that are already on the spool.

Fortunately for us,
the branch of topology known as braid theory
(itself a part of knot theory)
exactly addresses this situation.
That is (with the shade in its folded position),
we require the core of the band to be in \textit{braid position},
which we briefly summarize below.

\section*{Braids: A Brief Introduction}

\begin{figure}[ht]

 \centering
 \includegraphics[width=20mm]{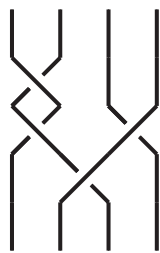}
 \caption{A four-strand braid.
 \label{braid-Figure}
 }
\end{figure}

Informally, a \dfn{braid} is a set of one or more strands of string
all starting at the same height as each other,
extending down without local minima or maxima but with twisting allowed,
and ending at the same height as each other.
Figure~\ref{braid-Figure} shows an example.
(For a more detailed introduction to braid theory, see Adams~\cite{adams}.)

\begin{figure}[ht]

 \centering
 \includegraphics[width=80mm]{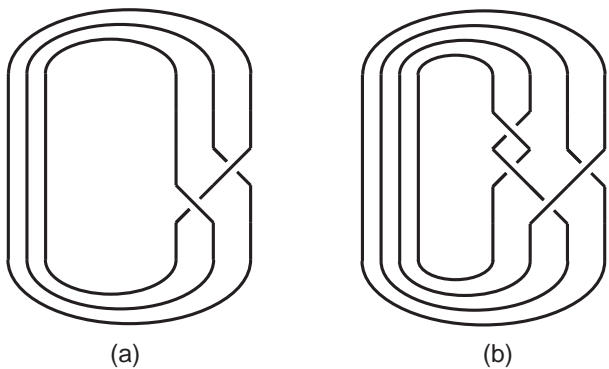}
 \caption{(a)~A one-component closed braid. (b)~A two-component closed braid.
 \label{two-closed-braids-Figure}
 }
\end{figure}

A \dfn{closed braid} is a braid that is ``closed up''
by connecting the top ends of the braid to its bottom ends
in a pairwise fashion.
Thus, the braid is actually made up of
one or more continuous loops of string,
each of which we call a \dfn{component}.
For example,
when the four-strand braid in Figure~\ref{braid-Figure} is closed up,
we end up with the two-component closed braid in
Figure~\ref{two-closed-braids-Figure}(b).
The strands are drawn off to the side to avoid introducing any new crossings into our diagram.

\begin{figure}[ht]

 \centering
 \includegraphics[width=45mm]{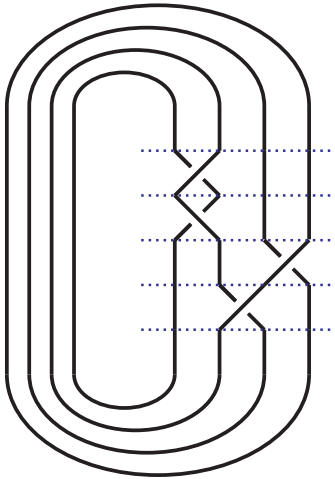}
 \caption{Braid position: only one crossing in each ``horizontal section''.
 \label{sliced-braid-Figure}
 }
\end{figure}

In any (closed) braid one can easily arrange that
whenever two strands cross each other,
no other pair of strands cross each other at that same height.
Thus, parallel horizontal lines can be drawn such that
between any two consecutive lines
there is exactly one crossing,
as in Figure~\ref{sliced-braid-Figure}.
When a closed braid is put in this configuration,
we say it is in \dfn{braid position}.
Note that we have now made mathematically precise
the idea of a ribbon wrapped on a spool.
Thus we have:
\begin{statement}
With the shade in its closed position, the core of the band
is a single-component closed braid in braid position.
\end{statement}

Now suppose our single-component closed braid was constructed
from a braid containing $m$ strands and $n$ crossings. It is
not difficult to see that in this case $m$ and $n$ must have
opposite parity. Figure~\ref{two-closed-braids-Figure}(a) shows
a simple example of this. We will not prove this fact here, but
we will give a brief justification for those who are familiar
with the symmetric groups: The $n$ crossings between the $m$
strands translate into a sequence of $n$ transpositions in
$S_m$. Having a \textit{single}-component closed braid
translates into saying that the product of these $n$
transpositions gives an $m$-cycle in $S_m$. And it is known
that in this case $m$ and $n$ must have opposite parity
(Lanski \cite[pages 140-141, Theorems~3.13, 3.14]{lanski} ). We now
have:

\begin{statement}
  With the shade in its closed position,
  if the band is coiled into $m$ loops and the core of the band
  (which is in braid position) has $n$ crossings,
  then $m$ and $n$ have opposite parity.
\end{statement}
Now it's time to consider the two boundary circles of the band.
This leads us naturally to another topic in topology,
namely links and linking number.

\section*{Linking Number: Another Brief Introduction}

\begin{figure}[ht]

 \centering
 \includegraphics[width=110mm]{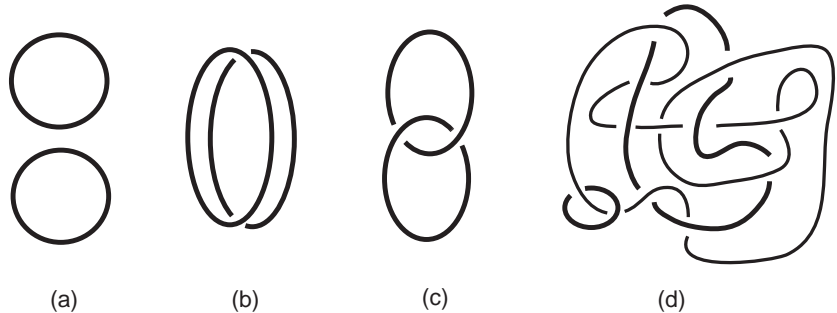}
 \caption{Three two-component links, and one three-component link.
 \label{link-examples-Figure}
 }
\end{figure}

Informally, a \dfn{link}
is a set of $n$ continuous loops
(again called components) in $3$-space.
Figure~\ref{link-examples-Figure} shows a few examples.
(When $n$ happens to be 1, we also call our link a knot.)
Note that a closed braid is
simply a specific presentation of a knot or link.
(For a more detailed introduction to links, see Adams~\cite{adams}.)

There is a very useful notion,
called the \dfn{linking number},
that applies to links with two components.
Roughly speaking,
it is the number of times
one of the components twists around the other.
We will require a precise definition,
which we develop below.

Given a diagram of a 2-component link,
we first orient each component,
i.e., we put a ``direction'' on each component of the link
(as if telling an ant which way
to travel).
The direction we choose for each component is arbitrary.

\begin{figure}[ht]

\vspace{3mm}
 \centering
 \includegraphics[width=40mm]{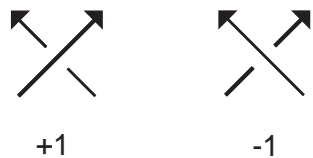}
 \caption{A right-handed ($+1$) and a left-handed ($-1$) crossing.
 \label{right-hand-rule-Figure}
 }
\end{figure}

Now, we focus only on crossings
between the two components
(rather than crossings between a component and itself).
Each such crossing consists of two oriented segments or ``arrows'', as shown in Figure~\ref{right-hand-rule-Figure}.
We assign $+1$ or $-1$ to each crossing as follows: If the head of the top arrow is to the right of the bottom arrow (when facing in the direction of the bottom arrow), we assign $+1$ to the crossing; otherwise we assign $-1$.

\begin{figure}[ht]

\vspace{3mm}
 \centering
 \includegraphics[width=40mm]{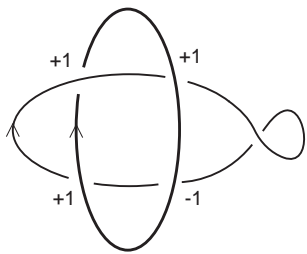}
 \caption{Linking number $= (1 + 1 + 1 - 1)/2 =1$.
 \label{linking-number-Figure}
 }
\end{figure}

To finish, we add up all the $+1$s and $-1$s
and then divide the sum by 2.
The result is called the linking number of the link
for the given orientations.
Figure~\ref{linking-number-Figure} shows an example of
an oriented link with linking number 1.

What will be important for us is the fact that
the absolute value of the linking number is a
\textit{topological invariant} of the link.
That is, no matter how the link is ``moved around,''
what diagram of it we look at,
and what orientation we pick for each component,
the absolute value of the linking number will not change.
The interested reader is referred to
Adams~\cite{adams} for a proof of this fact.
Now, with the shade in its open position,
the linking number of the boundary circles is clearly zero
(see Figures~\ref{open-band-Figure} and \ref{link-examples-Figure}(b)).
Since this linking number does not change
when the shade is folded up,
we have:
\begin{statement}
 With the shade in its closed position,
 the linking number of the boundary circles is zero.
\end{statement}

\section*{Analysis of the Magic Shade using Braid Position and Linking
Number}

Let's begin by recalling our first two Statements:
With the shade in its closed position,
we assume that the core of the band is
a single-component closed braid in braid position;
and if the band is coiled into $m$ loops
and the core of the band has $n$ crossings,
then $m$ and $n$ have opposite parity.

\begin{figure}[ht]

 \centering
 \includegraphics[width=100mm]{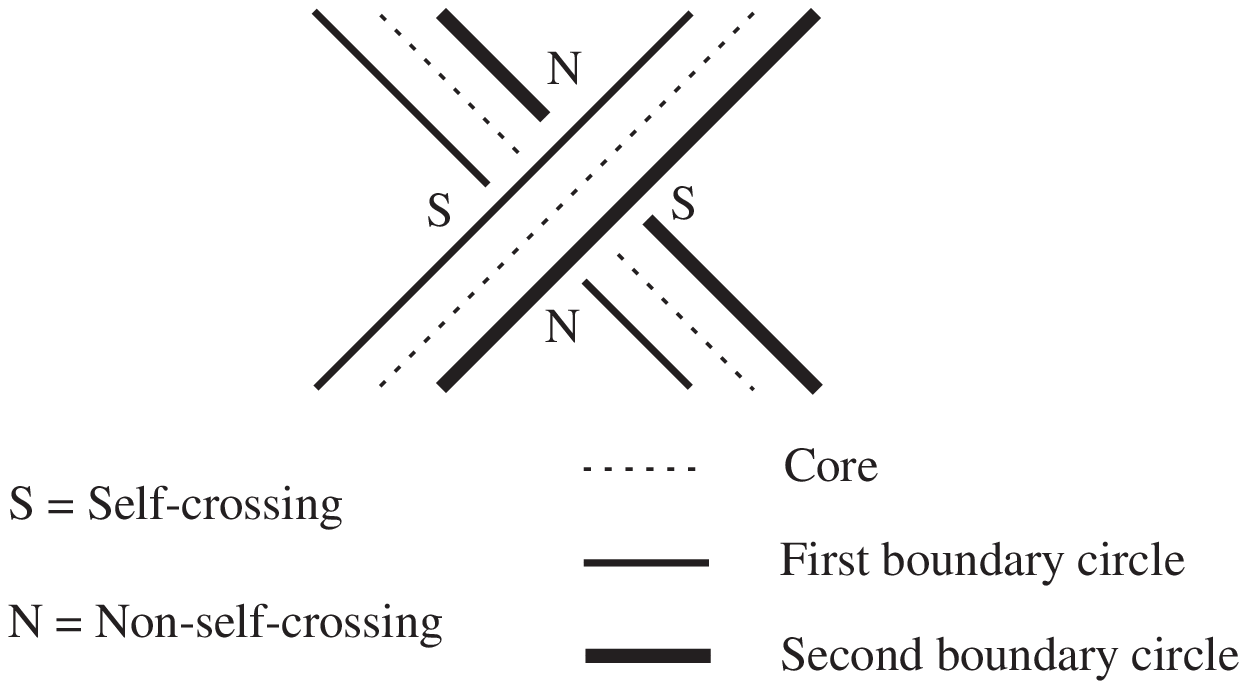}
 \caption{Two self- and two non-self-crossings of the
 boundary circles near one self-crossing of the core.
 \label{four-crossing-Figure}
 }
\end{figure}

Now draw in the two boundary circles, very close to the core.
Since the band has no twists in it,
the only crossings we will see
between the two boundary circles
are near the $n$ self-crossings of the core,
as in Figure~\ref{four-crossing-Figure}.
The core has now served its purpose;
let's take it out of the picture.
In place of each of the $n$ self-crossings of the core
we now see
four new crossings involving the boundary circles:
two self-crossings
(between a boundary circle and itself),
and two non-self-crossings
(between different boundary circles).
Since we're interested in the linking number of the boundary circles,
we ignore the former type of crossing
and concentrate only on the latter.
So, we have a total of $2n$ crossings,
arranged as $n$ pairs,
between the two boundary circles.

We compute the linking number of the two boundary circles as follows: First, we give each boundary circle an arbitrary orientation
and then assign $+1$s and $-1$s to each of the $2n$ crossings.
%But, according to our third Statement,
%the boundary circles have linking number zero.
%What does this tell us about $n$?
Next, recalling that the core of the band is in braid position, we note that each boundary circle is as well. Thus, looking again at Figure~\ref{four-crossing-Figure},
notice that for each boundary circle, its two strands in that figure
will be oriented either both downward or both upward.
It follows that in each of the $n$ pairs of crossings,
we have either two $+1$s or two $-1$s, never one of each.
Recall that the linking number is
the absolute value of \emph{half} the sum of all these
$+1$s and $-1$s;
hence each of these $n$ pairs of crossings contributes
exactly $1$ or $-1$, never $0$,
to the linking number.
Now, $+1$ and $-1$ are both odd, and if we add up $n$ odd numbers,
the result has the same parity as $n$.
But, by Statement~3, the linking number is zero;
so $n$ must be even.
Hence, by Statement~2,
$m$ must be odd!
That is:

\begin{theorem}
In its closed position,
the wire frame of a Magic Shade sunshade
must be coiled into an odd number of loops.
\end{theorem}

At this point,
if the reader has such a sunshade handy,
he or she should verify that the wire frame in fact coils into
\textit{three} loops when the shade is closed
--- thus answering Question~1 in the Introduction ---
and that indeed there are no twists in the wire.
To answer Question~2 in the Introduction,
the reader may also verify that
while it is possible to force the shade into two or four coils,
this will introduce twists in the wire,
and the shade will certainly not like being put in this position.
But what about larger odd numbers?
It turns out that with a little trial and error,
the shade can indeed be folded into five (quite small) loops!
The process for doing so actually involves a common topology ``trick'':
With the shade closed into three loops, make two of those loops very small.
Then close up the big loop
(as if it were an open shade)
into three loops,
while the two smaller ones ``go along for the ride''.
Following this idea, it is then clear that,
at least theoretically,
the shade could be folded into any odd number of loops.

\end{document}